\documentclass{amsart}[12pt]
\usepackage{xypic}

\usepackage{amsmath}
\usepackage{amstext}
\usepackage{amssymb}
\usepackage{amsthm}
\usepackage{amscd}

\usepackage{lscape}
\usepackage{longtable}

\usepackage{epsfig}
\input txdtools
\let\et=\etexdraw
\def\etexdraw{\drawbb\et}

\theoremstyle{plain}
\newtheorem{thm}{Theorem}[section]
\newtheorem{thm*}{Theorem}
\newtheorem{lem}[thm]{Lemma}
\newtheorem{prop}[thm]{Proposition}
\newtheorem{prop*}[thm*]{Proposition}
\newtheorem{cor}[thm]{Corollary}

\theoremstyle{definition}

\theoremstyle{remark}

\newtheorem{rmk}[thm]{Remark}

\DeclareMathOperator{\initial}{in} 
 \DeclareMathOperator{\lt}{lt}

\DeclareMathOperator{\height}{ht}
\DeclareMathOperator{\Ass}{Ass}

\DeclareMathOperator{\adj}{adj}

\begin{document}

\def\Ker {\operatorname{Ker}\nolimits}
\def\Im {\operatorname{Im}\nolimits}
\def\Image {\operatorname{Image}\nolimits}
\def\Syz {\operatorname{Syz}\nolimits}
\def\initial {\operatorname{in}\nolimits}
\def\Gin{\operatorname{Gin}\nolimits}
\def\Spec{\operatorname{Spec}\nolimits}
\def\D{\operatorname{D}\nolimits}
\def\V{\operatorname{V}\nolimits}
\def\H {\operatorname{H}\nolimits}
\def\E {\operatorname{E}\nolimits}
\def\V {\operatorname{V}\nolimits}
\def\nE {\operatorname{e}\nolimits}
\def\nV {\operatorname{v}\nolimits}
\def\C {\operatorname{\cal C}\nolimits}
\def\row {\operatorname{row}\nolimits}
\def\column {\operatorname{column}\nolimits}

\title
[On ideals of minors of matrices with indeterminate entries]
{On ideals of minors of matrices with indeterminate entries}
\author{Mordechai Katzman}
\address{Department of Pure Mathematics,
University of Sheffield, Hicks Building, Sheffield S3 7RH, United Kingdom\\
}
\email{M.Katzman@sheffield.ac.uk}


\date{\today}



\maketitle

\section{Introduction and notation}

This paper has two aims. The first is to study ideals of  minors of matrices whose entries are
among the variables of a polynomial ring. Specifically, we describe matrices whose ideals of  minors of a given size
are prime.
The ``generic'' case, where all the entries are distinct variables has been studied extensively
(cf. \cite{BC} and \cite{BV} for a thorough account.)
While some special cases, such as catalecticant matrices and other 1-generic matrices, have been studied
by other authors (e.g., \cite{E2}), the general case is not well understood.
The main result in the first part of this paper is
Theorem \ref{Some prime specialisations of determinantal ideals} which gives
sufficient conditions for the ideal of minors of a matrix to be prime.
This theorem is general enough to include interesting examples, such as the ideal of maximal minors of
catalecticant matrices and their generalisations discussed in the second part of the paper.

The second aim of this paper is to settle a
specific problem raised by David Eisenbud and Frank-Olaf Schreyer (cf.~\cite{ES}) on the primary decomposition
of an ideal of maximal minors. We solve this problem by applying  Theorem \ref{Some prime specialisations of determinantal ideals}
together with some ad-hoc techniques.

Throughout this paper $\mathbb{K}$ shall denote a field. For any matrix $M$ with entries in a ring and any
$t\geq 1$, $I_t(M)$ will denote the ideal generated by the $t\times t$ minors of $M$.
The results to be presented here rely on  well known properties of determinantal rings
which we summarise below:
\begin{thm}\label{recollections}
Let $X=(x_{ij})$ be the generic $m\times n$ matrix
and let $T=\mathbb{K}[x_{11}, \dots, x_{mn}]$.
\begin{enumerate}
    \item [(a)] $T/I_t(X)$ is a Cohen-Macaulay domain (cf. Theorems 1.10 and 6.7 in \cite{BC}),
    \item [(b)] $\dim T/I_t(X)=(m+n-t+1)(t-1)$ (cf. Theorem 1.10 in \cite{BC}).
\end{enumerate}
\end{thm}

\section{Some prime ideals of minors}

Throughout this section $X=\big(x_{ij}\big)$ will be a generic $m\times n$ matrix with $m\geq n$ and
$T$ will denote the polynomial ring over $\mathbb{K}$ whose variables are the entries of $X$.
We fix a $1\leq t\leq n$ and write $r=t-1$.
The aim of this section is to describe some prime ideals of minors of the image of $X$ under the identification
of some of the variables $x_{ij}$. We shall prove that these ideals are prime by embedding the appropriate quotient rings
into domains. As a first step we realise that determinantal varieties are rational:

\begin{prop}\label{determinantal rings are rational}
Let $W$ be a polynomial ring with variables $\{ y_{i j} \,|\, 1\leq i\leq m, 1\leq j\leq r\}$ and
$\{ z_{i j} \,|\, 1\leq i\leq r, 1\leq j\leq n-r\}$ over $\mathbb{K}$. Let $Y$ be a $m\times r$ matrix whose
$(i,j)$ entry is $y_{ij}$, let $Z$ be the following $r\times n$ matrix
$$
\left(
\begin{array} {llllllll}
1 & 0 & 0 & \dots & 0 & z_{1,1} & \dots & z_{1,n-r} \\
0 & 1 & 0 & \dots & 0 & z_{2,1} & \dots & z_{2,n-r} \\
  &   & \ddots &  &   & \vdots  &       & \\
0 & 0 & 0 & \dots & 1 & z_{r,1} & \dots & z_{r,n-r}
\end{array}
\right)
$$
and let $S$ be the sub-$\mathbb{K}$-algebra of $W$ generated by the entries of the product $Y Z$.
The map $T/I_t(X) \rightarrow S$ sending the image of $x_{ij}$ to the $(i,j)$ entry of $Y Z$ extends to an isomorphism
$f: T/I_t(X) \rightarrow S$ of $\mathbb{K}$-algebras.
\end{prop}
\begin{proof}
Notice that $Y Z$ has rank $r=t-1$, hence $f$ is well defined. Since $f$ is clearly surjective, we only need to show that
it is injective, and we show this by showing that $\dim S=\dim T/I_t(X)=(m+n-r) r$; we achieve this by
showing that $W$ and $S$ have the same fraction field.

Obviously, $y_{i j}/1$ is in the fraction field of $S$ for all $1\leq i\leq m, 1\leq j\leq r$;
we now
show that $z_{i j}/1$ is in the fraction field of $S$ for all $1\leq i\leq r$ and $1\leq j\leq n-r$.
Let $M$ be the $r\times r$ submatrix of $Y$ consisting of its first
$r$ rows and denote the classical adjoint of $M$ with $\adj M$.
Now the entries of $(\det \adj M)^{-1} (\adj M) M Z$ are in the fraction field of $S$
but $(\det \adj M)^{-1} (\adj M) M Z$ contains $\big(z_{i,j}\big)$ as a submatrix, and we are done.
\end{proof}

Let $X_r$ be the submatrix of $X$ consisting of its first $r$ columns.
We let $\mathcal{J}=\big\{ (i,j) \,|\, 1\leq i\leq m, 1\leq j\leq n \big\}$ and
$\mathcal{J}_r=\big\{ (i,j) \,|\, 1\leq i\leq m, 1\leq j\leq r \big\}$.

Given  a sequence
$\mathcal{S}=\big((\alpha_1,\beta_1), \dots, (\alpha_r,\beta_r)\big)$  elements of
$\mathcal{J}_r \times \mathcal{J}$ we define a sequence of
directed graphs $G_0, \dots, G_r$ whose vertex sets are $\mathcal{J}$ and whose edges given by
$$\E(G_j)= \Big\{ \overrightarrow{(a,b) (a,\ell)} \,|\, 1\leq a\leq m \, , 1\leq b\leq r \, , r+1\leq \ell\leq n \Big\}
\bigcup
\Big\{ \overrightarrow{\beta_1 \alpha_1 }, \dots, \overrightarrow{ \beta_j \alpha_j} \Big\}
.$$
We call the sequence $\mathcal{S}$ \emph{bad}, if for some $1\leq j\leq r$
there exists a directed path in $G_{j-1}$ starting at $\alpha_j$ and ending at $\beta_j$.

The motivation for this definition is as follows.

\begin{prop}\label{Proposition: bad sequences}
Use the notation of the previous proposition and write $A=Y Z$.
Let $\mathcal{S}=\big((\alpha_1,\beta_1), \dots, (\alpha_s,\beta_s)\big)$ be a set of elements
in $\mathcal{J}_r \times \mathcal{J}$ where
$\alpha_1, \dots, \alpha_s$ are distinct.

Define recursively a sequence of matrices $A^{(0)}, \dots, A^{(s)}$ as follows:
$A^{(0)}=A$ and, for all $0\leq i < s$, $A^{(i+1)}$ is obtained from $A^{(i)}$ by replacing
each occurrence of $y_{\alpha_{i+1}}=A^{(i)}_{\alpha_{i+1}}$ in $A^{(i)}$ with its $\beta_{i+1}$ entry.

If for some $\alpha\in\mathcal{J}_r$, $\beta\in\mathcal{J}$ and some $0\leq j \leq s$
$y_{\alpha}$ occurs in $A^{(j)}_\beta$ then there exists a directed path in $G_j$ starting
at $\alpha$ and ending at $\beta$.
\end{prop}
\begin{proof}
We proceed by induction on $j$; when $j=0$ notice that $y_{\alpha}$ occurs in $A_\beta$
for $\alpha\neq \beta$ if and only if
$\alpha=(a,b)$ and $\beta =(a,\ell)$ with $ 1\leq a\leq m$, $1\leq b\leq r$ and $r+1\leq \ell\leq n$.
Assume henceforth that $j>0$.

Pick a minimal $0\leq k \leq j$ such that $y_{\alpha}$ occurs in $A^{(k)}_\beta$.
If $k<j$, the induction hypothesis implies that
there exists a directed path in $G_k$ starting
at $\alpha$ and ending at $\beta$, and that path is also a directed path in $G_j$, and the theorem follows.

We assume now that $y_{\alpha}$ does not occur in $A^{(j-1)}_\beta$ but it does occur in
$A^{(j)}_\beta$. This implies that $y_{\alpha_j}$ occurs in $A^{(j-1)}_{\beta}$ and that
$y_{\alpha}$ occurs in $A^{(j-1)}_{\beta_j}$, and the induction hypothesis implies that
there exist in $G_{j-1}$ a directed path
$P_1$ starting at $\alpha$ and ending at ${\beta_j}$
and
a directed path
$P_2$ starting at $\alpha_j$ and ending at ${\beta}$.
Since both $P_1$ and $P_2$ are also directed paths in $G_j$,
and since $\overrightarrow{ \beta_j \alpha_j} \in \E(G_j)$
we have a directed path in $G_j$ from $\alpha$ to $\beta$ given by
the concatenation of the path $P_1$, followed by $\overrightarrow{ \beta_j \alpha_j}$ and $P_2$.
\end{proof}

\begin{thm}\label{Some prime specialisations of determinantal ideals}
Let $\mathcal{S}=\big((\alpha_1,\beta_1), \dots, (\alpha_s,\beta_s)\big)$ be a sequence of elements in
$\mathcal{J}_r \times \mathcal{J}$.
Let $L\subseteq T$ be the ideal minimally generated by
$x_{\alpha_1}-x_{\beta_1}, \dots, x_{\alpha_{s}}-x_{\beta_{s}}$.
Assume that
\begin{enumerate}
\item[(i)]
$\alpha_1, \dots, \alpha_s$ are distinct,
\item[(ii)]
$\mathcal{S}$ is not a bad sequence, and
\item[(iii)]
the image of
$I_{r}(X_r)$ in $ T/(L + I_t(X))$ has positive height.
\end{enumerate}
Then
$T/(L+I_t(X))$ is a Cohen-Macaulay domain of dimension $(m+n-t+1)(t-1)-s$.
Furthermore, $T/(L+I_t(X))$ is rational.
\end{thm}
\begin{proof}

Let $L^\prime\subseteq T$ be the ideal generated by
$x_{\alpha_1}-x_{\beta_1}, \dots, x_{\alpha_{s-1}}-x_{\beta_{s-1}}$.
Write
$$
U^\prime=\frac{T}{L^\prime+I_t(X)},\quad  U=\frac{T}{L+I_t(X)},$$
$$
W^\prime=\mathbb{K}[z_{i,j}] [y_\alpha \,|\, \alpha\in \mathcal{J}_r \setminus \{\alpha_1, \dots, \alpha_{s-1}\}],\quad
W=\mathbb{K}[z_{i,j}] [y_\alpha \,|\, \alpha\in \mathcal{J}_r \setminus \{\alpha_1, \dots, \alpha_{s}\}].
$$

Let $A$ be and $A^{(0)}, \dots, A^{(s)}$ be matrices as in Proposition \ref{Proposition: bad sequences}.
Write $B^\prime=(b^\prime_{ij})$ for $A^{(s-1)}$
and  $B=(b_{ij})$ for $A^{(s)}$.
Proposition \ref{Proposition: bad sequences} guarantees that
as we transform $B^\prime$ to $B$,
the expression  replacing $y_{\alpha_s}$
does not contain $y_{\alpha_s}$, and so the variable $y_{\alpha_s}$  does not occur in $B$.
Define the $\mathbb{K}$-algebra map $\rho: W^\prime \rightarrow W$
which sends $y_{\alpha_s}$ to the $\beta_s$ entry of $B^\prime$ and fixes all other variables.

Let $S^\prime$ be the sub-$\mathbb{K}$-algebra of $W^\prime$ generated by the entries of $B^\prime$ and
let $S$ be the sub-$\mathbb{K}$-algebra of $W$ generated by the entries of $B$.
We have the following commutative diagram
\begin{equation}\label{CD1}
\xymatrix{
U^\prime \ar@{->>}[r]^{f^\prime} \ar@{->>}[d]^{\phi} & S^\prime \ar@{^{(}->}[r]^{i^\prime} \ar[d]^{\psi} & W^\prime \ar@{->>}[d]^{\rho}\\
U \ar@{->>}[r]^{f} & S \ar@{^{(}->}[r]^{i} &  W
}
\end{equation}
where $f^\prime$ is the surjection which maps $x_{ij}$ to the $(i,j)$ entry of $B^\prime$
and $f$ is the surjection which maps $x_{ij}$ to the $(i,j)$ entry of $B$ (these are well defined
because the ranks of $B$ and $B^\prime$ are less than $t$)
and
where the restriction of $\rho$ to $S^\prime$
induces the map $\psi : S^\prime \rightarrow S$  which replaces any factor
$y_{\alpha_s}$ in a generator by the $\beta_s$ entry of $B^\prime$.
The map $\phi : U^\prime \rightarrow U=U^\prime/(x_{\alpha_s}-x_{\beta_s})$ is the quotient map.

To prove the theorem we
show that $f$ is an isomorphism and we do so by induction on $s$.
The case $s=0$ (i.e., $L=0$) is a restatement of
Proposition \ref{determinantal rings are rational}.
Assume now that $s>0$. Notice that the induction hypothesis implies that $U^\prime$ is a Cohen-Macaulay domain,
and so $U=U^\prime/(x_{\alpha_s}-x_{\beta_s})$ is also Cohen-Macaulay.

Write
$$
D^\prime=
\left[
\begin{array}{llll}
b^\prime_{11} & b^\prime_{12} & \dots & b^\prime_{1,r}\\
b^\prime_{21} & b^\prime_{22} & \dots & b^\prime_{2,r}\\
\vdots & \vdots &       & \vdots\\
b^\prime_{m1} & b^\prime_{m2} & \dots & b^\prime_{m,r}\\
\end{array}
\right], \quad
E^\prime=
\left[
\begin{array}{lll}
b^\prime_{1,r+1} & \dots & b^\prime_{1,n}\\
b^\prime_{2,r+1} & \dots & b^\prime_{2,n}\\
\vdots &        & \vdots\\
b^\prime_{m,r+1} & \dots & b^\prime_{m,n}\\
\end{array}
\right] .
$$
Notice that
$D^\prime \big( z_{i,j} \big)=E^\prime$
and that, if we localize at any non-zero $\delta^\prime\in I_r(D^\prime)$,
we can obtain each $z_{i,j}$ as a rational function of
the entries of $D^\prime$ and $E^\prime$.
Condition (i) implies that
$y_\alpha\in S^\prime$ for all $\alpha\in \mathcal{J}_r \setminus \{\alpha_1, \dots, \alpha_{s-1}\}$
and we deduce that $S^\prime_{\delta^\prime}=W^\prime_{\delta^\prime}$ for every such minor.

Since $f$ is clearly surjective, we conclude the proof by showing that $f$ is injective;
write $P=\ker f$.

The image of $I_r(X_r)$ in $T/(L+I_t(X))$ has positive height and since $T/(L+I_t(X))$ is Cohen-Macaulay
we can find a $d^\prime$ in the image of $I_r(X_r)$ in $T/(L^\prime+I_t(X))$ such that $d:=\phi(d^\prime)$
is  not a zero-divisor on  $T/(L+I_t(X))$. We now show that the localisation $P U_d$ vanishes and so
$P=0$.

Write $\delta^\prime=f^\prime(d^\prime)$ and
localise (\ref{CD1}) to obtain
\begin{equation*}
\xymatrix{
U^\prime_{d^\prime}  \ar@{->}[r]^{f^\prime_{d^\prime}}_\cong \ar@{->>}[d]^{\phi_{d^\prime}} & S^\prime_{\delta^\prime} \ar@{->}[r]^{i^\prime}_\cong & W^\prime_{\delta} \ar@{->>}[d]^{\rho_{\delta^\prime}}\\
U_d  &  &  W_{\delta}
}
\end{equation*}
Now the localisation of $i^\prime$ at $\delta^\prime$ is an isomorphism;
but
$\ker \rho_{\delta^\prime} : W^\prime_{d^\prime} \rightarrow W_{d}$ is generated by
$y_{\alpha_s}-B^\prime_{\beta_s}=f^\prime_{\delta^\prime}(x_{\alpha_s}-x_{\beta_s})$ so
$\ker f^\prime_{d^\prime} \circ i^\prime_{f^\prime(d^\prime)} \circ \rho_{\delta^\prime}$ is generated by
$x_{\alpha_s}-x_{\beta_s}$, i.e.,
$$\ker f^\prime_{d^\prime} \circ i^\prime_{f^\prime(d^\prime)} \circ \rho_{\delta^\prime}=\ker \phi_{d^\prime} .$$
Now $\phi_{d^\prime}^{-1} (PU_d) \subseteq \ker f^\prime_{d^\prime} \circ i^\prime_{f^\prime(d^\prime)} \circ \rho_{\delta^\prime}$
so $\phi_{d^\prime} \big( \phi_{d^\prime}^{-1} (P U_d) \big)=0 $ and we deduce that $P U_d=0$.
\end{proof}

Theorem \ref{Some prime specialisations of determinantal ideals} has numerous applications.
The following is an instance where conditions (i) and (ii) hold trivially.
\begin{cor}
Let $R_0$ be a polynomial ring. Fix $m\geq n$ and let $R$ be the polynomial ring over $R_0$ with variables $x_{ij}$ for
$1\leq i\leq m$, $t\leq j\leq n$ where $1\leq t \leq n$.
Let $A=(a_{ij})$ be a $m \times n$ matrix such that $a_{ij}$ are variables in $R_0$ for $j<t$ while
$a_{ij}=x_{ij}$ for $j\geq t$. Write $r=t-1$ and let $A_r$ be the submatrix of $A$ consisting of it first $r$ columns.
If the image of $I_{r}(A_r)$
in $R/I_t(A)$
has positive height then $R/I_t(A)$ is a Cohen-Macaulay domain.
\end{cor}

\begin{rmk}
Let $m=n=t=3$.
Notice that while condition (ii) fails for the matrix
$\left(
\begin{array}{lll}
a & X & b\\
X & c & X\\
d & e & f\\
\end{array}
\right)
$,
the conditions of the theorem apply to its transpose, and indeed its determinant is irreducible.
On the other hand, the matrix
$\left(
\begin{array}{lll}
a & X & b\\
X & c & X\\
d & X & f\\
\end{array}
\right)
$
cannot be manipulated to satisfy condition (ii) of the theorem whereas
its determinant is irreducible.
\end{rmk}

\begin{rmk}
Among the results describing when ideals of maximal minors of matrices are prime,
\cite{E2} contains the following one which is quite general.
\begin{quote}
{\bf Theorem:} Let $M$ be a $m\times n$ matrix of linear forms where $m\geq n$.
If $M$ is the image of a 1-generic matrix modulo $\leq n-2$ linear forms, then
$I_n(M)$ is prime.
\end{quote}
One may ask whether the hypothesis of Theorem \ref{Some prime specialisations of determinantal ideals}
imply those of this theorem in the case of maximal minors-- the answer is ``no'':
consider the matrix
$$M=\left(
\begin{array}{llll}
A& A& C &x_1\\
B& A& B &x_2 \\
A& C& B &x_3 \\
A& B& A &x_4 \\
\end{array}
\right) .
$$
A computation with using a Macaulay2 (\cite{GS}) script written by David Eisenbud shows that
$M$ is not the image of a 1-generic matrix modulo $\leq 2$ linear forms, while
the primality of $I_4(M)$ can be deduced from Theorem \ref{Some prime specialisations of determinantal ideals}.
\end{rmk}

\section{A primary decomposition}
Let $R=\mathbb{K}[a_1, a_2, a_3, a_4, a_5, b_1, b_2, b_3, b_4, b_5]$.
We let $$
M_5=\left(
\begin{array}{lllll}
a_1& a_2& a_3& a_4& a_5\\
b_1& b_2& b_3& b_4& b_5
\end{array}
\right),
M_4=\left(
\begin{array}{llll}
a_1& a_2& a_3& a_4\\
a_2& a_3& a_4& a_5\\
b_1& b_2& b_3& b_4\\
b_2& b_3& b_4& b_5
\end{array}
\right),
M_3=\left(
\begin{array}{lll}
a_1& a_2& a_3\\
a_2& a_3& a_4\\
a_3& a_4& a_5\\
b_1& b_2& b_3\\
b_2& b_3& b_4\\
b_3& b_4& b_5
\end{array}
\right).
$$
In a lecture in Overwolfach in April 2005, David Eisenbud conjectured that
that  $I_3(M_4)$ is radical with primary decomposition $I_2(M_5)\cap I_3(M_3)$
(see also section 4 of \cite{ES}).
In this section we show that this is indeed the case.

We shall need  the following elementary lemma:
\begin{lem}\label{lemma1}
Let $I$ be an homogeneous ideal in a polynomial ring $P=\mathbb{K}[y_1, \dots, y_n]$, and fix a term ordering in that ring.
Let $A\subset I$ be a finite set of homogeneous elements and let $\lt A$ be the set of leading terms of elements in $A$.
$\dim P/I \leq \dim P/P\lt A$.
\begin{proof}
Compute dimensions as one plus the degree of Hilbert polynomials,
recall that ideals and their initial ideals have identical Hilbert polynomials and
notice that $\lt A$ is contained in the initial ideal of $I$.
\end{proof}
\end{lem}


The following result could be obtained by proving that $M_3$ is 1-generic and applying Theorem 2.1 in \cite{E2}.
We give an alternative proof as an example of an application of
Theorem \ref{Some prime specialisations of determinantal ideals}.

\begin{cor}\label{I_3(M_3) is prime}
The ideal $I_3(M_3)$ is prime.
\end{cor}
\begin{proof}
We apply Theorem \ref{Some prime specialisations of determinantal ideals}
with $m=6$, $n=3$, $t=3$ the sequence of $x_\alpha$s taken to be
$x_{12}$,$ x_{13}$, $x_{22}$, $x_{23}$,
$x_{42}$, $x_{43}$, $x_{52}$, $x_{53}$
and sequence of $x_\beta$s taken to be
$x_{21}$, $x_{22}$, $x_{31}$, $x_{32}$,
$x_{51}$, $x_{52}$, $x_{61}$, $x_{62}$.
Now
$R/I_3(M_3)=T/(L+I_3(X))$.
It is not hard to verify that conditions (i) and (ii) of Theorem \ref{Some prime specialisations of determinantal ideals}
hold and it remains to verify that condition (iii) holds.

Define the matrix
\begin{eqnarray*}
D & = & \left(
\begin{array}{ll}
x & y \\
y & sx+ty \\
sx+ty & sy+t(sx+ty) \\
u & v \\
v & su+tv \\
su+tv & sv+t(su+tv)
\end{array}
\right)
\left(
\begin{array}{lll}
1 & 0 & s \\
0 & 1 & t
\end{array}
\right) \\
&=&
\left(
\begin{array}{lll}
x & y & sx+ty\\
y & sx+ty & sy+t(sx+ty)\\
sx+ty & sy+t(sx+ty) & s(sx+ty)+t(sy+t(sx+ty))\\
u & v & su+tv\\
v & su+tv & sv+t(su+tv)\\
su+tv & sv+t(su+tv) & s(su+tv)+t(sv+t(su+tv))
\end{array}
\right)
\end{eqnarray*}
and define $E$ to be the sub-$\mathbb{K}$-algebra of $\mathbb{K}[s,t,x,y,u,v]$
generated by the entries of $D$.
Define  $\phi: R \twoheadrightarrow E$ to be the surjection of $\mathbb{K}$-algebras which maps
the $(i,j)$ entry of $M_3$ to the $(i,j)$ entry of $D$ for all $1\leq i\leq 6, 1\leq j\leq 3$.
Notice that, if $\mathbf{c_1}, \mathbf{c_2}$ and $\mathbf{c_3}$ are the columns of $D$,
$\mathbf{c_3}=s \mathbf{c_1}+t \mathbf{c_2}$, and so
$I_3(M)\subseteq \Ker \phi$.
We now show that $\dim E=6$ by
showing that the fraction field $\mathbb{E}$ of $E$ is $\mathbb{K}(s,t,x,y,u,v)$, i.e., that
$s$ and $t$ are in  $\mathbb{E}$:
$$ \frac{[sx+ty]^2-y[sy+t(sx+ty)]}{x[sx+ty]-y^2}=\frac{s^2x^2+stxy-sy^2}{sx^2+txy-y^2}=s $$
and $t=([sx+ty]-sx)/y$.
Now since $I_3(M_3)\subseteq \ker \phi$ and $\dim E=6$, $\dim R/I_3(M_3) \geq 6$.

If we write
$$
N=
\left[
\begin{array}{ll}
a_1 & a_2 \\
a_2 & a_3 \\
a_3 & a_4 \\
b_1 & b_2\\
b_2 & b_3\\
b_3 & b_4
\end{array}
\right],
$$
since $I_3(M_3) \subseteq I_2(N)$ condition (iii) is equivalent to the statement
$\height I_3(M_3) < \height I_2(N)$.
We may choose a monomial order (say, reverse lexicographical) so that
the set of leading terms of the $2\times 2$ minors of $N$ contains $J=\{a_2^2 , a_3^2, a_4 b_1, b_2^2, b_3^2 \}$
and its not hard to see that $\dim R/RJ\leq 3$ hence Lemma \ref{lemma1} implies that
$\dim R/I_2(N) \leq 3$ and $\height I_2(N)\geq 8-3=5$.
But $\dim R/I_3(M_3)\geq 6$, so $\height I_3(M_3)\leq 10-6=4<\height I_2(N)$.
\end{proof}

\bigskip
Consider the permutations $\sigma, \tau$ of the variables of $R$
given by $\sigma(a_i)=a_{5-i+1}$, $\sigma(b_i)=b_{5-i+1}$ and
$\tau(a_i)=b_i$, $\tau(b_i)=a_i$ $(1\leq i\leq 5)$.
We shall use the fact that these can be extended to automorphisms of $R$ which fix
$I_2(M_5)$, $I_3(M_3)$ and $I_3(M_4)$.
We also denote henceforth $a_i b_j - a_j b_i$ with $\Delta_{ij}$.

\begin{prop}\label{delta{12} I_3(M_3) subseteq I_3(M_4)}
$\Delta_{12},\Delta_{23},\Delta_{13},\Delta_{45},\Delta_{34},\Delta_{35} \in \big( I_3(M_4) : I_3(M_3) \big)$.
\end{prop}
\begin{proof}

First notice that the generators of $I_3(M_3)$ not in $I_3(M_4)$ are
$$
d_1:=\left| \begin{array}{lll}
a_1 & a_2 & a_3\\
a_2 & a_3 & a_4\\
a_3 & a_4 & a_5
\end{array} \right|
,\quad
d_2:=\left| \begin{array}{lll}
a_1 & a_2 & a_3\\
a_2 & a_3 & a_4\\
b_3 & b_4 & b_5
\end{array} \right|
,\quad
d_3:=\left| \begin{array}{lll}
 a_1& a_2& a_3\\
 a_3& a_4& a_5\\
 b_2& b_3& b_4\\
\end{array} \right|
,\quad
d_4:=\left| \begin{array}{lll}
a_2 & a_3 & a_4\\
a_3 & a_4 & a_5\\
b_1 & b_2 & b_3
\end{array} \right|=\sigma\left( d_2 \right)
$$
and $\tau(d_1)$, $\tau(d_2)$,  $\tau(d_3)$, $\tau(d_4)$.
Indeed,
\begin{eqnarray*}
d_5:=\left| \begin{array}{lll}
a_1 & a_2 & a_3\\
a_3 & a_4 & a_5\\
b_1 & b_2 & b_3
\end{array} \right|
&=&
\left| \begin{array}{lll}
a_1 & a_2 & a_4\\
a_2 & a_3 & a_5\\
b_1 & b_2 & b_4
\end{array} \right|
-
\left| \begin{array}{lll}
a_1 & a_2 & a_3\\
a_2 & a_3 & a_4\\
b_2 & b_3 & b_4
\end{array} \right|
\in I_3(M_4)
,\\
d_6:=\left| \begin{array}{lll}
a_1 & a_2 & a_3\\
a_3 & a_4 & a_5\\
b_3 & b_4 & b_5
\end{array} \right|
&=&
\sigma\left( d_5 \right) \in I_3(M_4).\\
\end{eqnarray*}

Consider the relation
$$\left[
\begin{array}{lllll}
b_1& b_2& b_1& b_2& b_3\\
a_1& a_2& 0  & 0  & 0\\
0  &0   &a_1&a_2 &a_3\\
a_2&a_3 &a_2&a_3 &a_4\\
a_3&a_4 &a_3&a_4 &a_5\\
\end{array}\right]
\left[\begin{array}{l}
a_2\\ -a_1 \\ -a_2\\ a_1 \\ 0
\end{array}\right]
=0
$$
and call the $5\times 5$ matrix above $N_1$.
We expand $\det N_1$ using the first two columns to obtain
$$0=\det N_1 = \Delta_{12} d_1 -
\left| \begin{array}{ll}
a_1 & a_2 \\
a_2 & a_3
\end{array} \right| d_5 +
\left| \begin{array}{ll}
a_1 & a_2 \\
a_3 & a_4
\end{array} \right|
\left| \begin{array}{lll}
a_1 & a_2 & a_3\\
a_2 & a_3 & a_4\\
b_1 & b_2 & b_3
\end{array} \right| $$
hence $\Delta_{12} d_1 \in I_3(M_4)$.

Consider the relation
$$\left[
\begin{array}{lllll}
b_1& b_2& b_1& b_2& b_3\\
a_1& a_2& 0  & 0  & 0\\
0  &0   &a_1&a_2 &a_3\\
a_2&a_3 &a_2&a_3 &a_4\\
b_3&b_4 &b_3&b_4 &b_5\\
\end{array}\right]
\left[\begin{array}{l}
a_2\\ -a_1 \\ -a_2\\ a_1 \\ 0
\end{array}\right]
=0
$$
and call the $5\times 5$ matrix above $N_2$.
We expand $\det N_2$ using the first two columns to obtain
$$0=\det N_2 = \Delta_{12} d_2 -
\left| \begin{array}{ll}
a_1 & a_2 \\
a_2 & a_3
\end{array} \right| \tau(d_6) +
\left| \begin{array}{ll}
a_1 & a_2 \\
b_3 & b_4
\end{array} \right|
\left| \begin{array}{lll}
a_1 & a_2 & a_3\\
a_2 & a_3 & a_4\\
b_1 & b_2 & b_3
\end{array} \right| $$
hence $\Delta_{12} d_2 \in I_3(M_4)$.

Consider the relation
$$\left[
\begin{array}{lllll}
b_1& b_2& b_1& b_2& b_3\\
a_1& a_2& 0  & 0  & 0\\
0  &0   &a_1&a_2 &a_3\\
a_3&a_4 &a_3&a_4 &a_5\\
b_2&b_3 &b_2&b_3 &b_4\\
\end{array}\right]
\left[\begin{array}{l}
a_2\\ -a_1 \\ -a_2\\ a_1 \\ 0
\end{array}\right]
=0
$$
and call the $5\times 5$ matrix above $N_3$.
We expand $\det N_3$ using the first two columns to obtain
$$0=\det N_3 = \Delta_{12} d_3 -
\left| \begin{array}{ll}
a_1 & a_2 \\
a_3 & a_4
\end{array} \right|
\left| \begin{array}{lll}
a_1 & a_2 & a_3\\
b_1 & b_2 & b_3\\
b_2 & b_3 & b_4
\end{array} \right| +
\left| \begin{array}{ll}
a_1 & a_2 \\
b_2 & b_3
\end{array} \right| d_5 $$
hence $\Delta_{12} d_3 \in I_3(M_4)$.

Consider the relation
$$\left[
\begin{array}{lllll}
b_1& b_2& b_1& b_2& b_3\\
a_1& a_2& 0  & 0  & 0\\
0  &0   &a_1&a_2 &a_3\\
b_2&b_3 &b_2&b_3 &b_4\\
b_3&b_4 &b_3&b_4 &b_5\\
\end{array}\right]
\left[\begin{array}{l}
a_2\\ -a_1 \\ -a_2\\ a_1 \\ 0
\end{array}\right]
=0
$$
and call the $5\times 5$ matrix above $N_4$.
We expand $\det N_4$ using the first two columns to obtain
$$0=\det N_4 = \Delta_{12} \tau(d_4) -
\left| \begin{array}{ll}
a_1 & a_2 \\
b_2 & b_3
\end{array} \right|
\tau(d_5) +
\left| \begin{array}{ll}
a_1 & a_2 \\
b_3 & b_4
\end{array} \right|
\left| \begin{array}{lll}
a_1 & a_2 & a_3\\
b_1 & b_2 & b_3\\
b_2 & b_3 & b_4\\
\end{array} \right|$$
hence $\Delta_{12} \tau(d_4) \in I_3(M_4)$ and $\Delta_{12} d_4 \in I_3(M_4)$.

A similar argument employing the matrices
$$\left[
\begin{array}{lllll}
b_2& b_3& b_2& b_3& b_1\\
a_2& a_3& 0  & 0  & 0\\
0  &0   &a_2&a_3 &a_1\\
a_3&a_4 &a_3&a_4 &a_2\\
a_4&a_5 &a_4&a_5 &a_3\\
\end{array}\right]
,
\left[
\begin{array}{lllll}
b_2& b_3& b_2& b_3& b_1\\
a_2& a_3& 0  & 0  & 0\\
0  &0   &a_2&a_3 &a_1\\
a_3&a_4 &a_3&a_4 &a_2\\
b_4&b_5 &b_4&b_5 &b_3\\
\end{array}\right]
,
\left[
\begin{array}{lllll}
b_2& b_3& b_2& b_3& b_1\\
a_2& a_3& 0  & 0  & 0\\
0  &0   &a_2&a_3 &a_1\\
a_3&a_4 &a_3&a_4 &a_2\\
b_3&b_4 &b_3&b_4 &b_2\\
\end{array}\right]
,
\left[
\begin{array}{lllll}
b_2& b_3& b_2& b_3& b_1\\
a_2& a_3& 0  & 0  & 0\\
0  &0   &a_2&a_3 &a_1\\
a_3&a_4 &a_3&a_4 &a_2\\
b_4&b_5 &b_4&b_5 &b_3\\
\end{array}\right]$$
shows that $\Delta_{23} I_{3}(M_3) \subseteq I_3(M_4)$ and
a similar argument employing the matrices
$$\left[
\begin{array}{lllll}
b_1& b_3& b_1& b_3& b_2\\
a_1& a_3& 0  & 0  & 0\\
0  &0   &a_1&a_3 &a_2\\
a_2&a_4 &a_2&a_4 &a_3\\
a_3&a_5 &a_3&a_5 &a_4\\
\end{array}\right]
,
\left[
\begin{array}{lllll}
b_1& b_3& b_1& b_3& b_2\\
a_1& a_3& 0  & 0  & 0\\
0  &0   &a_1&a_3 &a_2\\
a_2&a_4 &a_2&a_4 &a_3\\
b_3&b_5 &b_3&b_5 &b_4\\
\end{array}\right]
,
\left[
\begin{array}{lllll}
b_1& b_3& b_1& b_3& b_2\\
a_1& a_3& 0  & 0  & 0\\
0  &0   &a_1&a_3 &a_2\\
a_2&a_4 &a_2&a_4 &a_3\\
b_2&b_4 &b_2&b_4 &b_3\\
\end{array}\right]
,
\left[
\begin{array}{lllll}
b_1& b_3& b_1& b_3& b_2\\
a_1& a_3& 0  & 0  & 0\\
0  &0   &a_1&a_3 &a_2\\
a_2&a_4 &a_2&a_4 &a_3\\
b_3&b_5 &b_3&b_5 &b_4\\
\end{array}\right]$$
shows that $\Delta_{13} I_{3}(M_3) \subseteq I_3(M_4)$.
We conclude the proof by noticing that
$\Delta_{45}=\sigma(\Delta_{12})$, $\Delta_{34}=\sigma(\Delta_{23})$ and $\Delta_{35}=\sigma(\Delta_{13})$.
\end{proof}

\begin{prop} \label{localisation at minors}
\begin{enumerate}
    \item[(a)] $I_3(M_4) \subseteq I_2(M_5) \cap I_3(M_3)$.
    \item[(b)] The ideal $I_3(M_4)$ is unmixed of height 4.
    \item[(c)] For any $P\in \Ass I_3(M_4) \setminus \{I_3(M_3)\}$,
    $I_2(M_5)\subseteq P$.
    \item[(d)] $\Ass I_3(M_4)= \{I_2(M_5), I_3(M_3)\}$.
\end{enumerate}
\end{prop}
\begin{proof}
(a)
Both inclusions $I_3(M_4) \subseteq I_2(M_5)$ and  $ I_3(M_4) \subseteq I_3(M_3)$
are easy to verify.

(b)
Recall from the proof of Corollary \ref{I_3(M_3) is prime} that $I_3(M_3)$ has height $4$.
From (a) we deduce that $\height I_3(M_4)\leq \height I_3(M_4)=4$,
and since $I_3(M_4) R_{\Delta_{12}} = I_3(M_3) R_{\Delta_{12}}$,
$\height I_3(M_4) \geq \height I_3(M_3) R_{\Delta_{12}}=4$.

Take $X$ to be the generic $4\times 4$ matrix as in Theorem \ref{recollections}.
That theorem tells us that $T/I_3(X)$ is Cohen-Macaulay and $12$-dimensional.
Also, if $J$ is the ideal of $T$ generated by the six elements
$x_{12}-x_{21}, x_{13}-x_{22}, x_{14}-x_{23}, x_{32}-x_{41}, x_{33}-x_{42}, x_{34}-x_{43}$,
we have $R/I_3(M_4)\cong P/(I_3(X)+J)$. But now
$$ 6 =  \dim R/I_3(M_4) = \dim T/(I_3(X)+J)$$
and so the six generators of $J$ form a system of parameters in the Cohen-Macaulay ring $T/I_3(X)$ and so
$R/I_3(M_4)\cong T/(I_3(X)+J)$ is Cohen-Macaulay, hence unmixed.

(c)
First notice that Proposition \ref{delta{12} I_3(M_3) subseteq I_3(M_4)}
implies that $\Delta_{12},\Delta_{23},\Delta_{13},\Delta_{45},\Delta_{34},\Delta_{35}\in P$.
We use the symmetry induced by $\sigma$, to reduce the problem to showing that
$P \supset \{ \Delta_{14}, \Delta_{24}, \Delta_{15} \}$.

Assume that $\Delta_{14}\notin P$:
$$ \left\{ \begin{array}{l}
a_3 \Delta_{14} + a_4 \Delta_{13} + a_1 \Delta_{34}=0\\
b_3 \Delta_{14} + b_4 \Delta_{13} + b_1 \Delta_{34}=0\\
\end{array}
\right.
\Rightarrow
 a_3, b_3 \in P$$
 and modulo $a_3, b_3$
$$
\left| \begin{array}{lll}
a_2 & a_3 & a_4\\
b_1 & b_2 & b_3\\
b_2 & b_3 & b_4
\end{array} \right| \equiv
\left| \begin{array}{lll}
a_2 & 0 & a_4\\
b_1 & b_2 & 0\\
b_2 & 0 & b_4
\end{array} \right|= b_2 \Delta_{24}
,\quad
\left| \begin{array}{lll}
a_1 & a_2 & a_3\\
a_2 & a_3 & a_4\\
b_2 & b_3 & b_4
\end{array} \right| \equiv
\left| \begin{array}{lll}
a_1 & a_2 &0\\
a_2 & 0 & a_4\\
b_2 & 0 & b_4
\end{array} \right|= a_2 \Delta_{24} .$$
If $\Delta_{24}\notin P$, we obtain $a_2, b_2\in P$;
now $P$ also contains the minor
$$\left| \begin{array}{lll}
a_1 & a_3 & a_4\\
a_2 & a_4 & a_5\\
b_1 & b_3 & b_4
\end{array} \right| \equiv a_4 \Delta_{14} \mod (a_2,a_3,b_2,b_3)$$
so $a_2,a_3,a_4,b_2,b_3\in P$, $\height P\geq 5$ and we obtain
a contradiction.
So now we assume that $\Delta_{24}\in P$, we deduce from
$$ \left\{ \begin{array}{l}
a_2 \Delta_{14} + a_4 \Delta_{12} + a_1 \Delta_{24}=0\\
b_2 \Delta_{14} + b_4 \Delta_{12} + b_1 \Delta_{24}=0\\
\end{array}
\right. $$
that $a_2, b_2\in P$ as well, and obtain, as before, a contradiction. We deduce that $\Delta_{14}\in P$.

A similar argument shows that $\Delta_{24}\in P$.

Assume now that $\Delta_{15}\notin P$:
$$ \left\{ \begin{array}{l}
a_3 \Delta_{15} + a_5 \Delta_{13} + a_1 \Delta_{35}=0\\
b_3 \Delta_{15} + b_5 \Delta_{13} + b_1 \Delta_{35}=0\\
a_4 \Delta_{15} + a_5 \Delta_{14} + a_1 \Delta_{45}=0\\
b_4 \Delta_{15} + b_5 \Delta_{14} + b_1 \Delta_{45}=0\\
\end{array}
\right.
\Rightarrow
a_3, b_3, a_4, b_4 \in P.$$
Modulo  $a_3, b_3, a_4, b_4$, the ideal $I_3(M_4)$ contains  the non-zero element $a_5 \Delta_{12}$ and,
if that were to happen, $\height I_3(M_4)>4$, a contradiction.

(d)
Both $I_3(M_3)$ and $I_2(M_5)$ are height-$4$ primes so they are minimal primes of $I_3(M_4)$ hence
$\{I_3(M_3), I_2(M_5)\} \subseteq \Ass I_3(M_4) $. Now (c) implies $ \Ass I_3(M_4)= \{I_3(M_3), I_2(M_5)\}$.

\end{proof}

\begin{prop}\label{d I_2(M_5) is in I3(M_4)}
Write $$\delta=
\left| \begin{array}{lll}
a_1 & a_2 & a_3\\
a_2 & a_3 & a_4\\
a_3 & a_4 & a_5
\end{array} \right| .$$
We have $\delta I_2(M_5) \subseteq I_3(M_4)$.
\end{prop}
\begin{proof}
Since
$$
\left| \begin{array}{llll}
a_1 & a_1 & a_2 & a_3\\
a_2 & a_2 & a_3 & a_4\\
a_3 & a_3 & a_4 & a_5\\
b_1 & b_1 & b_2 & b_3
\end{array} \right| =0
\Rightarrow
\left| \begin{array}{llll}
0 & a_1 & a_2 & a_3\\
a_2^2 -a_1 a_3 & a_2 & a_3 & a_4\\
a_2 a_3 -a_1 a_4& a_3 & a_4 & a_5\\
\Delta_{12} & b_1 & b_2 & b_3
\end{array} \right| =0
\Rightarrow$$
$$
-(a_2^2 -a_1 a_3)
\left| \begin{array}{lll}
a_1 & a_2 & a_3\\
a_3 & a_4 & a_5\\
b_1 & b_2 & b_3
\end{array} \right|
+
(a_2 a_3 -a_1 a_4)
\left| \begin{array}{lll}
a_1 & a_2 & a_3\\
a_2 & a_3 & a_4\\
b_1 & b_2 & b_3
\end{array} \right|
-\delta \Delta_{12}=0
$$
and
$$
\left| \begin{array}{lll}
a_1 & a_2 & a_3\\
a_3 & a_4 & a_5\\
b_1 & b_2 & b_3
\end{array} \right|=
\left| \begin{array}{lll}
a_1 & a_2 & a_4\\
a_2 & a_3 & a_5\\
b_1 & b_2 & b_4
\end{array} \right|
-
\left| \begin{array}{lll}
a_1 & a_2 & a_3\\
a_2 & a_3 & a_4\\
b_2 & b_3 & b_4
\end{array} \right|
\in I_3(M_4)
$$
we see that
$\delta \Delta_{12} \in I_3(M_4)$.
A similar argument shows that $\delta \Delta_{13}, \delta \Delta_{23} \in I_3(M_4)$.
By symmetry, i.e., an application of $\sigma$, we also obtain
$\delta \Delta_{45},\delta \Delta_{35}, \delta \Delta_{34} \in I_3(M_4)$.

Rather more mysteriously,
$$\delta \Delta_{14}=
\left| \begin{array}{lll}
a_1 & a_2 & a_3\\
a_2 & a_3 & a_4\\
b_1 & b_2 & b_3
\end{array} \right| (a_3 a_5 - a_4^2) +
\left| \begin{array}{lll}
a_1 & a_2 & a_3\\
a_2 & a_3 & a_4\\
b_2 & b_3 & b_4
\end{array} \right| (a_2 a_5 - a_3 a_4) +
\left| \begin{array}{lll}
a_1 & a_3 & a_4\\
a_2 & a_4 & a_5\\
b_1 & b_3 & b_4
\end{array} \right| (a_3^2 - a_2 a_4)$$
and hence also $\delta \sigma(\Delta_{14})=\delta \Delta_{25} \in I_3(M_4)$;
\begin{eqnarray*}
\delta \Delta_{15}&=&
-\left| \begin{array}{lll}
a_1 & a_2 & a_3\\
a_2 & a_3 & a_4\\
b_1 & b_2 & b_3
\end{array} \right| (a_3 a_5 - a_4^2) +
\left| \begin{array}{lll}
a_1 & a_2 & a_3\\
a_3 & a_4 & a_5\\
b_1 & b_2 & b_3
\end{array} \right| (a_2 a_5 - a_3 a_4)\\
&+&
\left| \begin{array}{lll}
a_1 & a_2 & a_4\\
a_2 & a_3 & a_5\\
b_2 & b_3 & b_5
\end{array} \right| (a_1 a_5 - a_3^2)+
\left| \begin{array}{lll}
a_1 & a_3 & a_4\\
a_2 & a_4 & a_5\\
b_1 & b_3 & b_4
\end{array} \right| (a_2 a_4-a_1 a_5)+
\left| \begin{array}{lll}
a_1 & a_3 & a_4\\
a_2 & a_4 & a_5\\
b_2 & b_4 & b_5
\end{array} \right| (a_2 a_3 -a_1 a_4);\\
\end{eqnarray*}
and
$$
\delta \Delta_{24}=
\left| \begin{array}{lll}
a_1 & a_2 & a_3\\
a_2 & a_3 & a_4\\
b_1 & b_2 & b_3
\end{array} \right| (a_3 a_5 - a_4^2) +
\left| \begin{array}{lll}
a_1 & a_2 & a_3\\
a_2 & a_3 & a_4\\
b_2 & b_3 & b_4
\end{array} \right| (a_2 a_5 - a_3 a_4)+
\left| \begin{array}{lll}
a_1 & a_3 & a_4\\
a_2 & a_4 & a_5\\
b_1 & b_3 & b_4
\end{array} \right| (a_3^2 - a_2 a_4).
$$

\end{proof}

\begin{thm}
The primary decomposition of $I_3(M_4)$ is given by
$I_2(M_5) \cap I_3(M_3)$.
\end{thm}
\begin{proof}
Proposition \ref{localisation at minors} shows that $\Ass I_3(M_4)=\{ I_2(M_5), I_3(M_3)\}$;
let $q_1 \cap q_2$ be the primary decomposition of  $I_3(M_4)$ where $q_1$ is associated to $I_2(M_5)$
and $q_2$ is associated to $I_3(M_3)$.

Let $\delta$ be as in Proposition \ref{d I_2(M_5) is in I3(M_4)};
$$q_1=I_3(M_4) R_{I_2(M_5)} \cap R \supseteq I_3(M_4) R_\delta \cap R \supseteq I_2(M_5)$$
hence $q_1=I_2(M_5)$.

An application of Proposition \ref{delta{12} I_3(M_3) subseteq I_3(M_4)} yields
$$q_2=I_3(M_4) R_{I_3(M_3)} \cap R \supseteq I_3(M_4) R_{\Delta_{12}} \cap R \supseteq I_3(M_3) $$
hence $q_2=I_3(M_3)$.

\end{proof}



\begin{thebibliography}{BV}

\bibitem{BC} W.~Bruns and A.~Conca.
\emph{Gr\"obner bases and determinantal ideals.}
Commutative algebra, singularities and computer algebra (Sinaia, 2002), pp.~9--66,
NATO Sci. Ser. II Math. Phys. Chem., 115.
Kluwer Acad. Publ., Dordrecht, 2003.

\bibitem{BV} W.~Bruns and U.~Vetter.
\emph{Determinantal rings.}
Lecture Notes in Mathematics, 1327.
Springer-Verlag, Berlin, 1988.

\bibitem{E1}
D.~Eisenbud.
\emph{Commutative algebra. With a view toward algebraic geometry.}
Graduate Texts in Mathematics, 150. Springer-Verlag, New York, 1995.

\bibitem{E2}
D.~Eisenbud.
\emph{Linear sections of determinantal varieties.}
American Journal of Mathematics {\bf 110} (1988), no. 3, pp.~541--575.



\bibitem{ES}
D.~Eisenbud and F-O.~Schreyer.
\emph{Relative Beilinson Monad and direct image for families of coherent sheaves.}
arXiv.math.AG/0506391

\bibitem{GS} D.~Grayson and M.~Stillman: Macaulay 2 --
a software system for algebraic geometry and commutative algebra,
available at {\tt http://www.math.uiuc.edu/Macaulay2}.



\end{thebibliography}
\end{document}